\documentclass[11pt]{article} % use larger type; default would be 10pt

\usepackage[utf8]{inputenc} % set input encoding (not needed with XeLaTeX)
\usepackage{url,color}
\usepackage[colorlinks=true,urlcolor=blue,linkcolor=blue,citecolor=magenta]{hyperref}
\usepackage{geometry} % to change the page dimensions
\usepackage{amsmath}
\usepackage{amssymb,amsthm,nicefrac,paralist}
\usepackage{graphicx} % support the \includegraphics command and options
\usepackage{mathtools} 
   
\usepackage{booktabs} % for much better looking tables
\usepackage{array} % for better arrays (eg matrices) in maths
\usepackage{paralist} % very flexible & customisable lists (eg. enumerate/itemize, etc.)
\usepackage{verbatim}  

\theoremstyle{plain}
\newtheorem{theorem}{Theorem}

\newtheorem{corollary}[theorem]{Corollary}
\newtheorem{proposition}[theorem]{Proposition}

\newtheorem*{theorem*}{Theorem}

\theoremstyle{definition}

\DeclarePairedDelimiter\floor{\lfloor}{\rfloor}
\newcommand\R{\mathbb{R}}
\newcommand\NC{\mathrm{NC}}
\title{Simple polytopes without small separators\footnote{The first author was funded by DFG through the \emph{Berlin Mathematical School}. Research by the second author was supported by the DFG Collaborative Research Center TRR~109 ``Discretization in Geometry and Dynamics.''}}
\author{Lauri Loiskekoski\\ 
Institut f\"ur Mathematik, FU Berlin\\Arnimallee 2\\14195 Berlin, Germany\\
\url{lauri.loiskekoski@fu-berlin.de}
\and  
G\"unter M. Ziegler\\
Institut f\"ur Mathematik, FU Berlin\\Arnimallee 2\\14195 Berlin, Germany\\
\url{ziegler@math.fu-berlin.de}}
\date{{\small October 2, 2015}\\[4mm]
\emph{For Gil Kalai on occasion of his sixtieth birthday}}

\begin{document}
\maketitle

\begin{abstract}
We show that by cutting off the vertices and then the edges of neighborly cubical polytopes,
one obtains simple $4$-dimensional polytopes with $n$ vertices 
such that all separators of the graph have size at least $\Omega(n/\log^{3/2}n)$. 
This disproves a conjecture by Kalai from 1991/2004.  
\end{abstract}

\section{Introduction}

The Lipton--Tarjan planar separator theorem from 1979 \cite{LiTa}
states that for any \emph{separation constant} $c$, with $0<c<\nicefrac12$, the vertex set of any
planar graph on $n$ vertices can be partitioned into three sets $A,B,C$ with $cn\le|A|\le|B|\le(1-c)n$
and $|C|=O(\sqrt n)$, such that $C$ separates $A$ from $B$, that is,
there are no edges between $A$ and~$B$.  Traditionally $c=\nicefrac13$ is used.
Miller, Thurston et al.\ \cite{MillerThurston-separators} in 1990/1991 provided a geometric proof for the planar separator theorem, 
combining the fact that every $3$-polytope has an edge tangent representation 
by the Koebe--Andreev--Thurston circle packing theorem with the center point theorem. 

Miller, Teng, Thurston \& Vavasis \cite{MillerTengThurstonVavasis} 
generalized the planar separator theorem to $d$ dimensions, that is,
to the intersection graphs of suitable ball packings in~$\R^d$. In view of this,
Kalai noted that there is no separator theorem for general $d$-polytopes, due to the existence of 
the cyclic polytopes, whose graph is complete for $d\ge4$ and thus has no separators.
However, he conjectured that the graphs of \emph{simple} $d$-polytopes 
cannot be good expanders, that is, they all should 
have small separators.
Specifically, in his 1991 paper on diameters and $f$-vector theory \cite[Conj.~12.1]{Ka2} 
(repeated in the 1997 first edition of~\cite{kalai04:_polyt})  
he postulated  that for every $d\ge3$ any simple $d$-polytope on $n$ vertices should have a separator of size
\[
  O\big( n^{1-\frac1{\floor{d/2}}}\big),
\]
which fails for $d=3$, while for $d=4$ it postulates separators of size $O(\sqrt n)$.
In the 2004 second edition of the \emph{Handbook} \cite[Conj.~20.2.12]{kalai04:_polyt} 
he corrected this to ask for separators of size
\[
  O\big( n^{1-\frac1{d-1}}\big),
\]
which for $d=3$ yields the planar separator theorem, and for $d=4$ 
conjectures the existence of separators of size $O(n^{2/3})$.

At that time Kalai also referred to \cite{MillerTengThurstonVavasis} for the claim that there are 
triangulations of~$S^3$ on $n$ tetrahedra
that cannot be even separated by $O(n/\log n)$ vertices. 
This is not stated in the paper \cite{MillerTengThurstonVavasis},
but it refers to a construction of Thurston who had described to his coauthors an embedding
of the cube-connected cycle graph in $\R^3$ as the dual graph of a configuration of tetrahedra.
Details about this construction seem to be lost (Gary Miller, personal communication 2015).
% If I remember correctly Thurston gave an embedding  of the cube-connected cycle graph in R^3 as linear tets.
% -- Gary Miller, Email April 8, 2015

In this note, we disprove Kalai's conjectures, and come close to confirming Thurston's claim.
Our construction uses the existence of \emph{neighborly cubical $4$-polytopes} $\NC_4(m)$,
first proved by Joswig \& Ziegler \cite{Z62}:
For each $m\ge4$ there is a $4$-dimensional polytope $\NC_4(m)$ 
whose graph is isomorphic to the graph $C_m$ of the $m$-cube
and whose facets are combinatorial $3$-cubes.

\begin{theorem}\label{thm:main}
    For any $m\ge4$ “cutting off the vertices, and then the original edges” from a neighborly
    cubical $4$-polytope $\NC_4(m)$ results in a simple $4$-dimensional polytope $\NC_4(m)''$ 
    with $n:=(6m-12)2^m$ vertices 
    whose graph has no separator of size less than 
    \[\Omega\Big(\frac n{\log^{3/2}n}\Big),\]
    while separators of size 
    \[O\Big(\frac n{\log n}\Big)\]
	exist.
\end{theorem}

To prove this, we only use the \emph{existence} of neighborly cubical
$4$-polytopes, and the fact that the $f$-vector of any such polytope is
\[ 
    f(\NC_4(m)) = (f_0,f_1,f_2,f_3) = 2^{m-2} (4, 2m, 3m-6, m-2), 
\]
but not a complete combinatorial description, as given in~\cite[Thm.~18]{Z62}.
Indeed, it was later established by Sanyal \& Ziegler \cite{Z102} that there are  
\emph{many} different combinatorial types, and Theorem~\ref{thm:main} and its proof
are valid for all of them.
It may still be that  
 for some specific neighborly cubical $4$-polytopes all separators in the resulting simple polytopes  
have size at least $\Omega(n/\log n)$. This would strongly confirm Thurston's claim.
On the other hand, all simple $4$-polytopes that are constructed according to Theorem \ref{thm:main}
have separators of size $O(n/\log n)$. However, we do not know whether such separators
exist for arbitrary simple $4$-polytopes on $n$ vertices.

This paper is structured as follows.
In Section~\ref{sec:construct} we describe the construction of $\NC_4(m)''$,  
compute its $f$-vector, establish that it is simple, and give a “coarse” description of the graph $C''_m:=G(\NC_4(m)'')$.
In Section~\ref{sec:separator} we show that the graph $C''_m$ has no small separators. 
This follows from elementary and well-known expansion properties of the cube graph $C_m$.
In Section~\ref{sec:upper_bound} we exhibit separators of size $O(n/\log n)$ in the graphs $C_m''$
derived from \emph{arbitrary} neighborly cubical $4$-polytopes. Finally, 
in Section~\ref{sec:comments} we extend all this to simple $d$-dimensional polytopes for $d\ge4$.

\section{Doubly truncated neighborly cubical polytopes}\label{sec:construct}
 
A \emph{neighborly cubical $d$-polytope} $\NC_d(m)$ is a $d$-dimensional convex polytope 
whose $k$-skeleton for $2k+2 \leq d$ is isomorphic to that of the $m$-cube. 
It is required to be \emph{cubical}, which means that all of its faces are combinatorial cubes. 
The existence of such polytopes was established by Joswig \& Ziegler \cite{Z62}.

For $4$-dimensional polytopes, the complete flag vector 
(that is, the extended $f$-vector of Bayer \& Billera \cite{BaBi}) 
is determined by the $f$-vector together
with the number $f_{03}$ of vertex-facet incidences.  

Let $m\ge4$. We start our constructions with a neighborly cubical $4$-polytope
$\NC_4(m)$ with the graph ($1$-skeleton) of the $m$-cube, so  
$f_0=2^m$ and $f_1=m2^{m-1}$. The rest of the flag vector is now
obtained from the Euler equation together with the fact that $\NC_4(m)$ is cubical:
Each facet has $6$ $2$-faces and $8$ vertices, which yields $6f_3=2f_{2}$
and $8f_3=f_{03}$. Thus we obtain
\begin{eqnarray*}
\textrm{flag} (\NC_4(m))&:=&(f_0, f_1, f_2, f_3; f_{03}) \\
&=& (2^m,m2^{m-1},3(m-2)2^{m-2},(m-2)2^{m-2}; 8(m-2)2^{m-2} )\\
&=& (4, 2m, 3m-6, m-2; 8m-16)\cdot2^{m-2}.
\end{eqnarray*} 

We generate the polytope $\NC_4(m)'$ from $\NC_4(m)$ by cutting off all of its vertices. 
The resulting polytope thus has
\begin{compactitem}[ $\bullet$ ]
    \item $(m-2)2^{m-2}$ facets that are $3$-cubes whose vertices have been cut off,
    with $f$-vector $(24,36,14)$,
    \item $2^m$ facets that are simplicial $3$-polytopes, each with $f$-vector $(m,3m-6,2m-4)$. 
\end{compactitem} 
The latter facets are the vertex figures of $\NC_4(m)$, which are
simplicial since the facets of $\NC_4(m)$, which are $3$-cubes, are simple. 
The resulting $4$-polytope has the following flag vector:
\begin{eqnarray*}
\textrm{flag} (\NC_4(m)')
&=& (4m, 14m-24, 11m-22,m+2; 28m-24)\cdot2^{m-2}.
\end{eqnarray*} 

We now generate $\NC_4(m)''$ from $\NC_4(m)'$ by cutting off the edges 
which come from edges in the original polytope $\NC_4(m)$ (but have been shortened in the transition to $\NC_4(m)'$). 
The resulting polytope has three types of facets:
\begin{compactitem}[ $\bullet$ ]
    \item $(m-2)2^{m-2}$ facets that are $3$-cubes whose vertices and subsequently
    the original cube edges have been cut off,
    with $f$-vector $(48,72,26)$,
    \item $2^m$ simple $3$-polytopes with $f$-vector 
    $(6m-12,9m-18,3m-4)$, which arise by cutting off the vertices of
    a simplicial $3$-polytope with $n$ vertices, and
    \item $m2^{m-1}$ prisms over polygons that may range between
    triangles and $(m-1)$-gons.  
\end{compactitem}
Again from the available information one can easily work out the flag vector, 
\begin{eqnarray*}
\textrm{flag} (\NC_4(m)'')
&=& (24m-48,48m-96,27m-46,3m+2; 28m-48)\cdot2^{m-2}.
\end{eqnarray*}
In particular, we can see from the $f$-vector that $f_1=2f_0$, so $\NC_4(m)''$
is simple.
Indeed, from \emph{any} $4$-polytope one gets a simple polytope by first cutting off the vertices, then 
the original edges. This may also be visualized in the dual picture: \emph{Any} $4$-polytope
may be made simplicial by first stacking onto the facets, and then
onto the ridges of the original polytope. After the first step,
the facets are pyramids over the original ridges. The second step 
corresponds to subdivisions of the pyramids in a point in the base,
which subdivides it into tetrahedra.

More generally, for $d$-polytopes we observe the following.

\begin{proposition}[see Ewald \& Shephard \cite{EwSh}]\label{prop:simplify}
    For $d\ge2$ and $0< k< d$ let $P$ be a $d$-polytope.
    Denote by $P^{(k)}$ the result of truncating the vertices, edges, etc.\
    up to the $(k-1)$-faces of $P$, in this order.
    Then the polytopes $P^{(d-2)}$ and $P^{(d-1)}$ are simple.
\end{proposition}
  
Indeed, in the dual picture stacking onto facets etc.\ down to edges, which yields $P^{(d-1)}$, 
corresponds to the barycentric subdivision of the boundary complex of the polytope. 
Subdividing the edges is unnecessary for our purpose, since these are already simplices.

\section{No small separators} \label{sec:separator}

Let $C_m$ be the graph of the $m$-cube, whose vertex set we identify with 
$\{0,1\}^m$. It has $2^m$ vertices and $m2^{m-1}$ edges. 
For any subset $S\subseteq V$ of the vertex set, 
its \emph{neighborhood} is defined as 
$N(S) := \{v \in V\setminus S : \{u,v\} \text{ is an edge for some }u\in S\}$.  
Harper solved the “discrete isoperimetric problem” 
in the $m$-cube in the sixties \cite{Harp}: For given
cardinality $|S|$, the cardinality of its neighborhood $|N(S)|$ is minimized by taking  
$S = \{v \in V, \sum_i^n v_i \leq d \}$ 
for some $d$ and taking the rest of vertices with coordinate sum $d+1$ in lexicographic order.
See Bollobás \cite{Bollobas-Combinatorics} or Harper \cite{Harper2004} for expositions. 
Thus optimal separators in the cube graph $C_m$ are obtained by
taking level sets, of size $\binom{m}k$.
Here the usual asymptotics for binomial coefficients
(as given by the central limit theorem, see e.g.\ \cite[Sect.~6.4]{Spencer:Asymptopia})
tell us that all the separators in the cube graph~$C_m$ have cardinality at least
$\Omega(2^m/\sqrt{m})$, where the implied constant depends on
the separation constant~$c$.
 
The graph $C'_m$ is obtained from the cube graph $C_m$ by replacing each node by a 
maximal planar graph on $m$ vertices and $3m-6$ edges. (Note that this description does
not specify the graph $C'_m$ completely.)
In the transition from $C'_m$ to $C''_m$, the $2^m$ planar graphs grow into cubic ($3$-regular) graphs
on $2(3m-6)=6m-12$ vertices each, which we call the \emph{clusters} of~$C_m''$. 
Each of the $m2^{m-1}$ edges between two vertices of~$C_m$ resp.\ between the
maximal planar graphs in $C'_m$ gives rise to
a number of edges (at least~$3$, at most $m-1$) between the corresponding two clusters in $C_m''$. 
While the cube graph~$C_m$ has $m2^{m-1}$ edges,
the modified graph $C_m''$ has $(6m-12)2^{m-1}$ edges between clusters.
Thus in the transition from $C_m$ to $C_m''$, the cube graph edges are replaced by less than $6$ edges on average.
 
$C_m''$ is a $4$-regular graph on $n=(6m-12)2^m$ vertices.
Consider an arbitrary separator of $C_m''$, consisting of  two disjoint sets of vertices $A$ and $B$ 
with $cn \le |A|\le |B|\le(1-c)n$ and a set~$C$ that contains the remaining vertices. 
From this we can generate a separator for the cube graph $C_m$
by labeling its vertices with $a$ or $b$ if the corresponding cluster in $C_m''$ has vertices only in $A$ or only in $B$, respectively. 
The remaining vertices will be labeled with $c$. 
There cannot be any neighboring vertices in $C_m$ labeled by $a$ and $b$, since this would imply 
neighboring $A$ and $B$ clusters in $C_m''$. 
The set of vertices of $C_m$ labeled $a$ has size at most $(1-c)n$, and the same is true for the 
set of vertices labeled $b$. Thus, for any fixed $c'<c$, unless the set of vertices labeled $c$ has linear size
(and thus we are done), both the sets of vertices with labels $a$ and $b$ have size at least $c'n$,
and thus we have constructed a separator for $C_m$.
By the isoperimetric inequality for vertex neighborhoods, there must be $\Omega(2^m/\sqrt{m})$ vertices labeled with $c$ and hence at least as many vertices in the separator for $C_m''$.
 
Thus all separators for the graph $C''_m$ of $\NC_4(m)''$  have size at least 
\[
\Omega\Big(\frac{2^m}{\sqrt m}\Big) \ =\ \Omega\Big(\frac{n}{\log^{3/2} n}\Big).
\]  

\section{Small separators} \label{sec:upper_bound}

Here we argue that for \emph{any} neighborly cubical $4$-polytope $\NC_4(m)$, the
derived simple $4$-polytope $\NC_4(m)''$ on $n=(6m-12)2^m$ vertices
has a separator of size 
\[
O(2^m) \ =\ O\Big(\frac{n}{\log n}\Big).
\]
Indeed, with respect to the identification of the vertex set of $C_m$ with $\{0,1\}^m$,
choose a random coordinate (“edge direction”), and 
divide the vertices of $C_m$ into two sets by whether the corresponding vertex label is 
$0$ or $1$. This corresponds to cutting the $m$-cube into two $(m-1)$-cubes,
with $2^{m-1}$ edges between them.

This cutting also divides the vertex set of $C_m''$ into two equal halves, containing
$n/2=(3m-6)2^m$ vertices each. In $C_m''$, there is an average of less than 6 edges between
adjacent clusters. For a random coordinate direction, the expected number
of edges between the two equal halves of $C_m''$ is less than $6\cdot 2^{m-1}=3\cdot 2^m$.
Thus by choosing a suitable coordinate, and removing one end vertex of each edge of 
$C_m''$ in the corresponding direction, we obtain a separator of
size less than $3\cdot 2^m$. 

\section{More generally} \label{sec:comments} 

We can extend the result of Theorem~\ref{thm:main} to dimensions $d>4$ by taking the product of 
$\NC_4(m)''$ and the standard $(d-4)$-cube. For a fixed dimension $d$ this gives 
a sequence of polytopes with $2^{d-4}$ times as many vertices. We can find a 
separator in this graph by taking a product of a separator in $\NC_4(m)''$ and 
the standard $(d-4)$-cube, so these polytopes are at least as easy to separate as 
$\NC_4(m)''$. On the  other hand the graph of this polytope again has a cube-like structure,
with $2^m$ clusters that are products of a cubic planar graph on $6m-12$ vertices
with the fixed graph $C_{d-4}$. Again we need to remove at least 
\[
\Omega\Big(\frac{2^m}{\sqrt m}\Big) \ =\ \Omega\Big(\frac{n}{\log^{3/2} n}\Big).
\] 
vertices to separate it. 

\begin{corollary}
    For each $d\ge4$ there is a sequence of simple $d$-polytopes whose graphs
    (on $n$ vertices) have no separators that are smaller than
    \[
    \Omega\Big(\frac{n}{\log^{3/2} n}\Big),
    \]
    but which have separators of size
    \[
    O\Big(\frac{n}{\log  n}\Big).
    \]    
\end{corollary}

Alternatively, one could try to start with neighborly cubical $d$-polytopes $\NC_d(m)$
and to “simplify” them by Proposition~\ref{prop:simplify}.
The resulting simple $d$-polytopes have graphs that are again similar to those of $m$-cubes,
where however the clusters have a size of the order of $\Theta(m^{\lfloor d/2\rfloor-1})$
for fixed $d$, and thus we get $n=\Theta(m^{\lfloor d/2\rfloor-1}2^m)$ vertices in total,
and thus separators of size  
\[
O(2^m)  \ =\ O\Big(\frac{n}{\log^{\lfloor d/2\rfloor-1}n}\Big).
\]
So the product construction sketched above is better for $d>5$. 

\subsubsection*{Acknowledgements.}
We are grateful to Hao Chen and Arnau Padrol for useful discussions.

\bibliographystyle{siam}%{plain}%{unsrt}   % this means that the order of references
			    % is dtermined by the order in which the
			    % \cite and \nocite commands appear
\begin{small}

%\bibliography{separators} 

\end{small}

\end{document}